\documentclass[a4paper,12pt]{article} 

\usepackage{amssymb}
\usepackage{epsfig}
\usepackage[english]{babel}
\usepackage{color}
\newcommand{\equ}[1]{(\ref{#1})}
\def\mmp{{{\cal M}_{\lambda,\Lambda}^+}}
\def\mmpb{{{\cal M}_{\bar\lambda,\Lambda}^+}}
\def\mmpn{{{\cal M}_{\lambda_n,\Lambda}^+}}
\def\mmm{{{\cal M}_{\lambda,\Lambda}^-}}
\def\mmpj{{{\cal M}_{\lambda_*,\Lambda}^+}}
\def\mmmj{{{\cal M}_{\lambda_*,\Lambda}^-}}
\def\mmpm{{{\cal M}_{\lambda,\Lambda}^\pm}}

\def\NN{{I\!\! N}}
\def\RR{I\!\!R}

\def\beq{\begin{equation}}
\def\eeq{\end{equation}}
\newcommand{\be}{\begin{equation}}
\newcommand{\ee}{\end{equation}}

\def\lm{L^{2}(\Omega)}
\def\li{C(\Omega)}
\def\lu{C([0, 1])}

\def\Box{\hfill\framebox(0.25,0.25){}}

\newtheorem{teo}{Theorem}[section]

\newtheorem{prop}{Proposition}[section]

\newtheorem{Remark}{Remark}[section]
\newcommand{\re}[1]{(\ref{#1})}
\newtheorem{lem}{Lemma}[section]
\newtheorem{cor}{Corollary}[section] 

\begin{document}

\vspace*{1cm} {\large\centerline{\bf Nonlinear Eigenvalues and Bifurcation 
Problems}
\centerline{\bf for Pucci's Operator }}

\bigskip
\bigskip
\centerline{ by}
\bigskip
\centerline{ J\'er\^ome BUSCA$^{*}$, Maria J. ESTEBAN$^{*}$}
\bigskip
\centerline{and}
\bigskip
\centerline{ Alexander QUAAS$^{**}$}
\bigskip
{\small
\centerline{$^{*}$Ceremade UMR CNRS 75341, Universit\'e Paris IX- Dauphine }
\centerline{ 75775 Paris Cedex 16, France.}
\centerline{\tt busca, esteban@ceremade.dauphine.fr} 

\bigskip
\centerline{$^{**}$ Departamento de Ingenier\'{\i}a  Matem\'atica,}
\centerline{ and CMM, UMR2071, CNRS-UChile}
\centerline{ Universidad de Chile, Casilla 170 Correo 3,}
\centerline{ Santiago, Chile.}
\centerline{\tt quaas11@dim.uchile.cl }
} 

\bigskip
\bigskip
\bigskip 

\centerline{\bf Abstract} 

\bigskip 

In this paper we extend existing results concerning generalized
eigenvalues of Pucci's extremal operators. In the radial case, we
also give a complete description of their spectrum, together with
an equivalent of Rabinowitz's Global Bifurcation Theorem. This
allows us to solve equations involving Pucci's operators. 

\bigskip
\bigskip
This work was partially supported by ECOS grant No. C02E08. The
third author is supported by CONICYT Becas de Postgrado.
\footnote{{\bf AMS Subject Classification.} Primary 35J60, 35P05,
35B32. Secondary 34B18.\\}

\setcounter{equation}{0}
\section{Introduction} 

If the solvability of fully nonlinear elliptic equations of the
form \beq\label{general_problem} F(x,u,Du,D^2u)=0 \eeq has been
extensively investigated for {\em coercive} uniformly elliptic
operators $F$, comparatively little is known when the assumption
on coercivity (that is, monotonicity in $u$) is dropped. 
In this paper, we want to
focus on the model problem \beq\label{model_problem} \left\{
\begin{array}{rcll}
 -\mmp(D^2u)&=& f(u) \quad\mbox{in} \quad \Omega ,\\
\\
      u&=&0  \quad  \mbox{on}\quad  \partial \Omega,\\
\end{array}\right.\eeq
(resp. $\mmm$) where $\Omega$ is a bounded regular domain, and
$\mmpm$ are the extremal Pucci's operators \cite{CPucci1} with
parameters $0<\lambda\leq \Lambda$  defined by
\begin{eqnarray*}
\mmp(M) = \Lambda \sum\limits_{e_i>0}e_i
+\lambda\sum\limits_{e_i<0}e_i
\end{eqnarray*}
and
\begin{eqnarray*}
\mmm(M) = \lambda \sum\limits_{e_i>0}e_i
+\Lambda\sum\limits_{e_i<0}e_i,
\end{eqnarray*}
for any symmetric $N\times N$ matrix $M$; here $e_i = e_i(M),
i=1,...,N,$ denote the eigenvalues of $M$. We intend to study
(\ref{model_problem}) as a bifurcation problem from the trivial
solution. Since $\mmpm$ are homogeneous of degree one, it is
natural to investigate the associated ``eigenvalue problem''
\beq\label{EV} \left\{
\begin{array}{rcll}-\mmp(D^2u)&=& \mu u \quad\mbox{in} \quad \Omega ,\\
\\
      u&=&0  \quad  \mbox{on}\quad  \partial \Omega.\\
\end{array}\right.
\eeq 
(resp. $\mmm$)  Pucci's extremal operators appear in the context of stochastic
control when the diffusion coefficient is a control variable, see
the book of Bensoussan and J.L. Lions \cite{ben} or the papers of
P.L. Lions \cite{lionspI}, \cite{lionspII}, \cite{lionspIII} for the relation beetween a general Hamilton-Jacobi-Bellman
and stochastic control.
They also provide natural extremal equations in the sense that if $F$
in (\ref{general_problem}) is uniformly elliptic, with ellipticity
constants $\lambda$, $\Lambda$, and depends only on the Hessian
$D^2u$, then 

\beq\label{ineg_extremale} \mmm(M) \leq F(M) \leq \mmp(M) \eeq for
any symmetric matrix $M$.

When $\lambda=\Lambda=1$, $\mmpm$ coincide with the Laplace
operator, so that (\ref{model_problem}) reads 

\beq\label{semilinear}
\left\{
\begin{array}{rcll}- \Delta u&=& f(u) \quad\mbox{in} \quad \Omega ,\\
\\
      u&=&0  \quad  \mbox{on}\quad  \partial \Omega,\\
\end{array}\right.\eeq 

whereas \re{EV} simply reduces to 

\beq\label{EVL} \left\{
\begin{array}{rcll}-\Delta u &=& \mu u \quad\mbox{in} \quad \Omega ,\\
\\
      u&=&0  \quad  \mbox{on}\quad  \partial \Omega.\\
\end{array}\right.
\eeq 

It is a very well known fact that there exists a sequence of
solutions
$$\{(\mu_n, \varphi_n)\}_{n \geq 1}$$
to \re{EVL} such that: 

i) the eigenvalues $\{\mu_n\}_{n \geq 1}$ are real, with $\mu_n>0$
and $\mu_n \to \infty$ as $n \to \infty $; 

ii) the set of all eigenfunctions $\{ \varphi_n\}_{n \geq 1}$ is a
basis of $\lm$. 

Building on these eigenvalues, the classical Rabinowitz
bifurcation theory \cite{rabinowitz2}, \cite{rabinowitz} then provides a comprehensive
answer to the existence of solutions of (\ref{semilinear}). 

When $\lambda<\Lambda$, problems (\ref{model_problem})-(\ref{EV})
are fully nonlinear. It is our purpose to investigate to which
extent the results about the Laplace operator can be generalized
to this context. A few partial results in this direction have been
established in the recent years and will be recalled shortly.
However, they are all concerned with the first eigenvalue and
special nonlinearities $f$. We provide here  a bifurcation result for general nonlinearities
 from the first two
``half-eigenvalues''  in general bounded domain. And in the radial case a
complete description of the spectrum and the bifurcation branches 
for a general nonlinearity from any point in the spectrum. 

\bigskip
Let us mention that besides the fact that
(\ref{model_problem})-(\ref{EV}) appears to be a favorable case
from which one might hope to address general problems like
(\ref{general_problem}), there are other reasons why one should be
interested in Pucci's extremal operators or, more generally, in
Hamilton-Jacobi-Bellman operators, which are envelopes of linear
operators. As a matter of fact, the problem under study has some
relation to the Fu\v c\'{\i}k spectrum. To explain this, let $u$
be a solution of the following problem
$$-\Delta u = \mu u^+-\alpha \mu u^-,$$
where $\alpha$ is a fixed positive number. One easily checks that
if $\alpha \geq 1$, then $u$ satisfies
$$\max\{ -\Delta u \, ,\,\frac{-1}{\alpha}\Delta u\}= \mu u,$$
whereas if $\alpha \leq 1$, $u$ satisfies
$$\min\{ -\Delta u \, ,\,\frac{-1}{\alpha}\Delta u\}= \mu u.$$
These relations mean that the Fu\v c\'{\i}k spectrum can be seen
as the spectrum of the maximum or minimum of two linear operators,
whereas (\ref{model_problem})-(\ref{EV}) deal with an infinite
family of operators.

We observe that understanding all the ``spectrum'' of the above
problem is essentially the same as determining the Fu\v c\'{\i}k
spectrum, which in dimension $N \geq 2$ is still largely an open
question, for which only partial results are known and, in general, they refer to a region near the usual spectrum,
(that is for $\alpha$ near 1). For a further discussion of
this topic, we refer the interested reader to the works of de
Figueiredo and Gossez \cite{Gossed}, H. Berestycki \cite{beres},
E.N. Dancer \cite{dancer}, S. Fu\v c\'{\i}k \cite{fucik}, P. Dr\'abek  
\cite{dra}, T.Gallouet and O. Kavian \cite{kavian2}, M. Schechter
\cite{schechter} and the references therein.

\bigskip
Our first result deals with the existence and characterizations of
the two first ``half-eigenvalues''. Some parts of it are already known (see below), but some are new.

\begin{prop}\label{fet}
Let $\Omega$ be a regular domain.
There exist two positive constants $\mu^+_1,\,\mu^-_1$, that we call first
half-eigenvalues such that: 

\bigskip
i) There exist two functions $\varphi^+_1,\, \varphi^-_1 \in C^2(\Omega)
\cap C(\bar\Omega)$
such that $(\mu^+_1,\,\varphi^+_1) $,  $(\mu^-_1,\,\varphi^-_1) $
are solutions to  \re{EV} and  $\varphi^+_1>0, \varphi^-_1<0 \quad
\mbox{in}\quad \Omega$.
Moreover, these two half-eigenvalues are simple, that is, all positive
solutions to \re{EV}
are of the form $(\mu^+_1,\, \alpha\varphi^+_1)$, with $\alpha>0$.
The same holds for the negative solution. 

\bigskip
ii) The two first half-eigenvalues satisfy
$$\mu^+_1=\displaystyle\inf_{A \in {\cal A}}\mu_1(A), \quad \mu^-_1=\sup_{A
\in {\cal A}}\mu_1(A),$$
where ${\cal A}$ is the set of all symmetric measurable matrices such that
$0<\lambda I \leq A(x) \leq \Lambda I$
and  $\mu_1(A)$ is the principal eigenvalue of the nondivergent second order linear elliptic
operator associated to $A$.

\bigskip 
iii) The two half-eigenvalues have the following characterization
$$\mu^+_1=\displaystyle\sup_{u>0
}\mathop{\mathrm{essinf}}_{\Omega}(-\frac{\mmp (D^2 u)}{u}), \quad
\mu^-_1=\sup_{u<0 }\mathop{\mathrm{essinf}}_{\Omega}(-\frac{\mmp
(D^2 u)}{u}).$$ The supremum is taken over all functions $u \in
W^{2,N}_{\mbox{loc}}(\Omega) \cap C(\bar\Omega)$. 

\bigskip
iv) The first half-eigenvalues can be also characterized by
$$\mu^+_1=\sup\{\mu\,|\,\mbox{ there exists }\phi>0\,\mbox{in } \Omega\quad
\mbox{satisfying}\quad \mmp (D^2 \phi)+\mu \phi \leq 0\}$$
 $$\mu^-_1=\sup\{\mu\,|\,\mbox{ there exists }\phi<0\,\mbox{in } \Omega\quad
\mbox{satisfying}\quad \mmp (D^2 \phi)+\mu \phi \geq 0\}$$
\end{prop}

\begin{Remark} \label{rem1}
Here and in the sequel, unless otherwise stated, it is
implicitly understood that any solution (resp. sub--,
super-solution) satisfies the corresponding equation (inequation)
pointwise a.e. This is the framework of {\em strong solutions}
\cite{GT}. 
\end{Remark}

The above existence result, that is part i) of Proposition
\ref{fet}, can been easily proved using an adaptation, for convex (or concave) operators, 
of  Krein-Rutman's Theorem in positive cones (see \cite{felmerquaas2} in the radial symmetric case 
and see \cite{quaas} in regular bounded domain). 

This existence result, has been proved recently in the case of general positive homogeneous
fully nonlinear elliptic operator, see the paper of Rouy
\cite{rouy}. The method used there is due to P.L. Lions who proved
the result i) of Proposition \ref{fet} for the Bellman operator (see
\cite{lions1}) and for the Monge-Amp\`ere operator (see
\cite{lions2}). Moreover, the definition of $\mu^+_1$ there translates  in our
case as:
\beq \label{rem12}
\mu^+_1=\sup\{ \mu \,|\, \mu \in {\cal I}\},
\eeq
where
$${\cal I}=\{ \mu \,|\, \exists \,\phi>0 \;\;\mbox{ s. t. } \;
\phi=0\;\;\mbox{on}\;\; \partial\Omega,\quad \mmp (D^2 \phi)+\mu \phi=-1 \;\; \mbox{ in }\;\; \Omega\}. $$ 

Properties ii) of  Proposition \ref{fet} can be generalized to any fully
nonlinear elliptic operator $F$ that is positively homogeneous of
degree one, with ellipticity constants $\lambda,\,\Lambda $. This
follows by the proof of ii) and (\ref{ineg_extremale}). 
These properties  were established by C. Pucci in
\cite{CPucci}, for other kind of extremal operators, see the
comments in Section 2.

The characterization of the form iii) and iv) for the first
eigenvalue, were introduced by Berestycki, Nirenberg and Varadhan
for second order linear elliptic operators (see \cite{varadan}). 

From the characterization iv) it follows that
$$\mu^+_1(\Omega) \leq \mu^+_1(\Omega')\quad\mbox{and}\quad \mu^-_1(\Omega) 
\leq \mu^-_1(\Omega')\quad\mbox{if}\quad \Omega' \subset \Omega.$$ 

For the two first half-eigenvalues, many other properties
will be deduced from the previous Poposition (See section 2). For example, 
whenever $\lambda\not=\Lambda$, we have $\,\mu_1^+ < \mu_1^-$, since
$\,\mu_1^+\leq \lambda \,\mu_1(-\Delta) \leq \Lambda\, \mu_1(-\Delta)\leq \mu_1^-$.

Another interesting and useful property is the following maximum principle.

\begin{teo}\label{fet1}
The next two maximum principles hold: 

\bigskip
a) Let $u \in W^{2,N}_{\mbox{loc}}(\Omega) \cap C(\bar\Omega)$ satisfy
\begin{eqnarray}\label{Mp1}
\begin{array}{rcll}\mmp(D^2u)+\mu u& \geq & 0 \quad\mbox{in} \quad \Omega
,\\
      u& \leq &0  \quad  \mbox{on}\quad  \partial \Omega.\\
\end{array}
\end{eqnarray}
If $\mu < \mu^+_1$, then $u \leq 0$ in $\Omega$. 

\bigskip
b) Let $u \in W^{2,N}_{\mbox{loc}}(\Omega) \cap C(\bar\Omega)$ satisfy
\begin{eqnarray}\label{Mp2}
\begin{array}{rcll}\mmp(D^2u)+\mu u& \leq & 0 \quad\mbox{in} \quad \Omega
,\\
      u& \geq &0  \quad  \mbox{on}\quad  \partial \Omega.\\
\end{array}
\end{eqnarray}
If $\mu < \mu^-_1$, then $u \geq 0$ in $\Omega$. 
\end{teo} 

\begin{Remark}
 These maximum principles are still valid for
continuous solutions in $\bar\Omega$ that satisfy the respective
inequalities in the viscosity sense. 
\end{Remark}

Results like Proposition \ref{fet} and  Theorem \ref{fet1}  can be obtained for $\mmm$ and can be
deduced just by noting that $\mmp(-M)=-\mmm(M)$, for any symmetric
matrix $M$. 

Next, we want to look  at the higher eigenvalues of  Pucci's extremal operators.
For that purpose we restrict ourselves to the radial case. 
In this case we have a precise description of the whole  ``spectrum'' and 
we expect that the result below will shed some light on the general case.
 More precisely, we have the following theorem.

\begin{teo}\label{teovpr}Let $\Omega=B_1$.
The set of all the scalars $\mu$ such that \re{EV} admits a nontrivial
radial solution,
consists of two unbounded increasing sequences
$$0<\mu^+_1<\mu^+_2<\cdot\cdot\cdot<\mu^+_k<\cdot\cdot\cdot,$$
$$0<\mu^-_1<\mu^-_2<\cdot\cdot\cdot<\mu^-_k<\cdot\cdot\cdot.$$
Moreover, the set of radial solutions of \re{EV} for $\mu=\mu^+_k$ is positively
spanned
by a function $\varphi^+_k$, which is positive at the origin and has exactly
k-1 zeros in $(0,1)$,
all these zeros being simple.
The same holds for $\mu=\mu^-_k$, but considering $\varphi^-_k$ negative at the origin.
\end{teo} 

Finally, we want to adress our original motivation, that is, 
we want to prove existence results for an equation of the type \re{model_problem}.
For this purpose we consider the nonlinear bifurcation problem associated with the extremal Pucci's operator, that is
\begin{eqnarray}\label{EVB}
\begin{array}{rcll}-\mmp(D^2u)&=& \mu u +f(u, \mu)\quad\mbox{in} \quad
\Omega ,\\
      u&=&0  \quad  \mbox{on}\quad  \partial \Omega,\\
\end{array}
\end{eqnarray}
where $f$ is continuous, $f(s, \mu)=o(|s|)$ near $ s=0$, uniformly for
$\mu \in \RR$ and $\Omega$ is a general bounded domain.
Concerning this problem we have the following theorem

\begin{teo}\label{bteo}

The pair $(\mu^+_1,0)$ (resp. $(\mu^-_1, 0)$) is a bifurcation point of positive (resp. negative) solutions to
\re{EVB}.
Moreover, the set of nontrivial solutions of \re{EVB} whose closure
contains $(\mu^+_1, 0)$ (resp. $(\mu^-_1, 0)$), is either unbounded
or contains a pair $(\bar\mu, 0)$ for some $\bar\mu$, eigenvalue of \re{EV}
with $\bar\mu\not=\mu^+_1$ (resp. $\bar\mu\not=\mu^-_1$).

\end{teo} 

Notice that a similar theorem can be proved in the case of $\mmm$.
The difference with Theorem \ref{bteo} is that $(\mu^+_1, 0)$ will be a
bifurcation
point for the negative solutions
and $(\mu^-_1, 0)$ will be a bifurcation point for the positive solutions.

\begin{Remark} Figure 1 allows to visualize the above result in which the 
bifurcation generates
only ``half-branches'': $u>0$ or $u<0$ in $\Omega$. 
\end{Remark}

\begin{figure}
\begin{center} 

\input{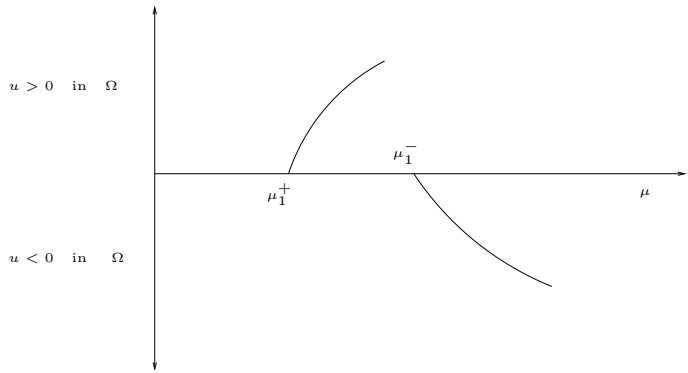}
\caption{ Bifurcation diagram for the first half-eigenvalues in a general 
bounded domain.}

\end{center}
\end{figure}

For the Laplacian the result is well known, see \cite{rabinowitz2}, \cite{rabinowitz},\cite{rabinowitz1}. In this case the ``half-branches'' become connected.
Therefore, we observe a symmetry breaking phenomena when $\lambda < \Lambda$.

For the $p$-Laplacian the result is known, in the general case, 
see the paper of del Pino and Man\'asevich \cite{Manuel}. See also the paper of del Pino, 
Elgueta and Man\'asevich \cite{elgueta}, for the case $N=1$. In this case the branches are also connected.
The proof of these results uses an invariance under homotopy with respect
to $p$ for the Leray-Schauder degree. 
In our proof of Theorem \ref{bteo} we use instead  homotopy invariance with respect to $\lambda$ (the ellipticity constant),
having to deal with a delicate region in which the degree is equal to zero.

A bifurcation result in the particular case $f(u, \mu)=-\mu |u|^{p-1}u$ can be 
found in the paper by P.L. Lions for the Bellman equation \cite{lions1}. 
For the problem

$$-\mmp(D^2u)=\mu g(x,u) \quad\mbox{in}\quad\Omega,\quad u=0
\quad\mbox{on}\quad\partial\Omega $$
with the following assumption on $g$:

\noindent
i) $u \to g(x,u)$ is nondecreasing and $g(x,0)=0,$ 

\noindent
ii) $u \to \frac{g(x,u)}{u}$ decreasing, and 

\noindent
iii)$$\displaystyle \lim_{u \to 0} \frac{g(x,u)}{u}=1, \quad \lim_{u \to
\infty} \frac{g(x,u)}{u}=0$$
a similar result was proved by E. Rouy  \cite{rouy}. 

In \cite{lions1} and \cite{rouy} the assumptions made were used in a
crucial way to construct sub and super solutions.
By contrast, we use a Leray-Schauder degree argument which allows us to treat
general nonlinearities. 

\bigskip
\noindent 
Other kind of existence results for positive solution of \re{model_problem}, can be
found in \cite{felmerquaas}, \cite{felmerquaas1}, \cite{felmerquaas2} and 
\cite{quaas}.

In the radially symmetric case we obtain a more complete result.
\begin{teo}\label{teobr} Let $\,\Omega=B_1$. For each $k \in \NN$, $k \geq 1$ there are two connected components
$S^\pm_k$ of nontrivial solutions to \re{EVB},
whose closures contains $(\mu^\pm_k,0)$.
Moreover, $S^\pm_k$ are unbounded and $(\mu, u) \in S^\pm_k$ implies that $u$
possesses exactly $k-1$ zeros in $(0,1)$.
\end{teo}
\begin{Remark}
1) $S^+_k$ (resp. $S^-_k$) denotes the set of solutions that are positive (resp. negative) at 
the origin.

\bigskip
\noindent
2) Figure 2, allows to visualize the above result in which the 
bifurcation generates
only ``half-branches'': $u(0)>0$ or $u(0)<0$. 
\end{Remark}
\begin{figure}
\begin{center} 

\input{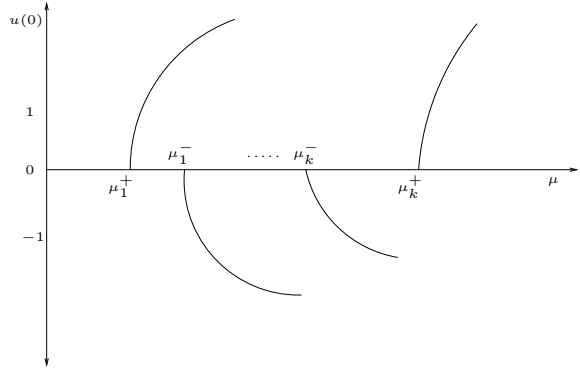}
\caption{ Bifurcation diagram in the radially symmetric case (note that $\mu^+_1 \leq \mu^-_1$, 
but for $k \geq 2$ the ordering between $\mu^+_k$ and $\mu^-_k$ is not known).}

\end{center}
\end{figure} 

For the Laplacian this result is well known. In this case, for all 
$k\geq 1$, $\mu^+_k=\mu^-_k$ and the `` half-branches '' now become connected. 

Our proof is based on the invariance of the Leray-Schauder degree under homotopy. It also uses some non-existence results.

The paper is organized as follows. In Section 2 we study the problem in a general regular bounded domain,
there we prove Theorem 1.1 and Theorem 1.3.
In section 3 we study the radial symmetric case, and we prove Theorem 1.2 and Theorem 1.4.

\section{First ``Eigenvalues'' in a General Domain and Nonlinear
Bifurcation}
We shall need the following version of Hopf's boundary lemma. 

\begin{lem}\label{lemh} Let $\Omega$ be a regular domain and let $u \in
W^{2,N}_{\mbox{loc}}(\Omega) \cap C(\bar\Omega)$ be a
non-negative solution to \beq \label{equa11} \mmm(D^2 u)+ \mu u
\leq 0 \quad\mbox{in}\quad \Omega,\quad u=0\quad \mbox{on}\quad
\partial \Omega, \eeq with $\mu \in \RR$. Then $u(x)>0$ for all
$x\in \Omega$. Moreover,
$$\displaystyle\limsup_{x \to x_0}\frac{u(x_0)-u(x)}{|x-x_0|}<0,$$
where $x_0 \in \partial \Omega$ and the limit is non-tangential,
that is, taken over the set of $x$ for which the angle between
$x-x_0$ and the outer normal at $x_0$ is less than $\pi/2-\delta$
for some fixed $\delta>0$.
\end{lem}
\noindent
\begin{Remark}\label{rhopf}
1) For a general strong maximum principle for degenerate convex elliptic
operators, see
the paper of M. Bardi, F. Da Lio \cite{bardi}. 

\bigskip
\noindent
2) This Lemma holds also for $u \in  C(\bar\Omega)$ that
satisfies the equation \re{equa11} in the viscosity sense. 

\bigskip
\noindent
3) A solution of \re{EV} in a regular domain is
necessarily $C^{2,\alpha}$ up to the boundary, see \cite{S}. Thus,
if $u$ is a positive solution to \re{EV}, then we have
$\frac{\partial u}{\partial \nu}<0$ on $\partial \Omega$ (resp.
$u<0$, $\frac{\partial u}{\partial \nu}>0$) if $\nu$ denotes the
outer normal.
\end{Remark} 

\noindent {\bf Proof}. We use the classical Hopf barrier function,
see for instance Lemma 3.4 in \cite{GT}. The rest of the proof
follows the lines of this lemma by using the weak maximum
principle, of P.L. Lions \cite{lions3} for solutions in
$W^{2,N}_{\mbox{loc}}(\Omega)$. $\Box$ 

\bigskip
\noindent Now we are in position to prove Proposition \ref{fet}. 

\smallskip
\noindent {\bf Proof of Proposition \ref{fet}.} i) The existence and simplicity
follow by using a Krein-Rutman's Theorem in positive
cones, see \cite{quaas}. For alternative methods see \cite{lions1} and
\cite{rouy}. Notice that by above Remark part 3), the two first 
half-eigenfunction are
$C^{2,\alpha}(\bar\Omega)$.

\bigskip
\noindent
ii) First notice that for a fixed function $v \in 
W^{2,N}_{\mbox{loc}}(\Omega)$ there exists
a symmetric
measurable matrix
$A(x)\in {\cal A} $, such that
$$\mmp(D^2 v)=L_A v,$$
where $L_A$ is the second order elliptic operator associated to $A$ , see 
\cite{CPucci1}. 

That is
$$L_A=\sum A_{i,j}(x)\partial_{i,j}.$$
Since $\varphi^+_1 \in C^2(\Omega)$, $\mu^+_1 \geq \inf_{A \in
{\cal A}}\mu_1(A)$. Suppose now for contradiction that $$\mu^+_1 >
\inf_{A \in {\cal A}}\mu_1(A).$$
Hence, there exists
$\bar A \in {\cal A}$ such that $\mu^+_1>\mu_1(\bar A)$.
The corresponding eigenfunction $u_1$ satisfies
$$-L_{\bar A} u_1=\mu_1(\bar A)u_1.$$
Moreover, $u_1 \in C^{2, \alpha}(\bar\Omega)$ and
$\frac{\partial u_1}{\partial \nu}<0$ on $\partial \Omega$. Same holds for 
$\varphi_1^+$. Thus there exists $K>0$
such that $u_1<K\varphi_1^+$. Notice that $u_1$ is a
sub-solution and $\varepsilon\varphi_1^+$ is a super-solution to
$$ \mmp(D^2 v)+\mu v=0 \quad\mbox{in}\quad \Omega.$$
Hence using Perron's method we find a positive solution to \re{EV},
which is in contradiction with part i).
Perron's method in this setting can be found for example in \cite{lions1}. 

\bigskip
\noindent
iii) We only need to prove that for any positive function $\phi \in
W^{2,N}_{\mbox{loc}}(\Omega) \cap C(\bar\Omega)$ we have 
$$\mu^+_1 \geq \inf_\Omega \frac{-\mmp(D^2 \phi)}{\phi}$$
Suppose the contrary, then there exists a positive function
$u \in   W^{2,N}_{\mbox{loc}}(\Omega) \cap C(\bar\Omega)$ and $\delta>0$
such that
$$\mu^+_1 +\delta < \inf_\Omega \frac{-\mmp(D^2 u)}{u}$$
So, $u$ satisfies
$$\mmp(D^2 u)+(\mu^+_1 +\delta)u \leq 0\quad\mbox{in}\quad \Omega .$$ 

On the other hand, $\varphi^+_1$ satisfies
$$\mmp(D^2 \varphi^+_1 )+(\mu^+_1 +\delta)\varphi^+_1 \geq
0\quad\mbox{in}\quad \Omega.$$
Using Lemma \ref{lemh} and Remark \ref{rhopf}(3), we can find $\varepsilon>0
$ such that
$$\varepsilon \varphi_1 \leq u \quad \mbox{in}\quad \Omega $$ 

Then, using Perron's method we find a positive solution to the problem
$$\mmp(D^2 v)+(\mu^+_1 +\delta)v=0\quad\mbox{in}\quad \Omega,$$
contradicting the uniqueness of the positive solution to \re{EV}, part i). 

\bigskip
\noindent
iv) follows directly from iii).$\Box$ 

\bigskip 

{\bf Proof of Theorem \ref{fet1}.}

 Let $u$ by a solution \re{Mp1} and $\hat{A} \in {\cal A}$
by such that $L_{\hat{A}}(u)=\mmp(D^2u)$.
Define $\hat{L}(v):=L_{\hat{A}}(v)+\mu v$. Using ii) of Proposition \ref{fet},, $ \mu_1(\hat{A})\geq
\mu^+_1>\mu$.
Then, clearly, the first eigenvalue of  $\hat{L}$ is positive.
So the maximum principle holds for $\hat{L}$ see \cite{varadan}. That is if
$v$ satisfies
\begin{eqnarray}\label{IP}
\begin{array}{rcll}\hat{L}(v) & \geq & 0 \quad\mbox{in} \quad \Omega ,\\
      v& \leq &0  \quad  \mbox{on}\quad  \partial \Omega.\\
\end{array}
\end{eqnarray}
implies $v \leq 0$. Since $u$ satisfies \re{IP}, $u \leq 0$. 

The same kind of argument can be used in case b).$\Box $

Now we will recall the following compactness results
for the Pucci's extremal operator, whose proof can be found for instance in
\cite{cab}. 

\begin{prop}\label{pc} Let $\{f_n\}_{n>0} \subset C(\Omega)$ be a
bounded sequence and\\
$\{u_n\}_{n>0} \subset C(\bar\Omega)\cap W^{2,N}_{loc}(\Omega)$ be a
sequence of solutions to
\beq\label{eqcom}
\mmp(D^2u_n)\geq f_n \quad\mbox{and}\quad \mmm(D^2u_n)\leq
f_n\quad\mbox{in}\quad \Omega,
\eeq
$u_n=0$ on $\partial\Omega$.
Then, there exists $u \in C(\bar\Omega)$ such that, up to a subsequence,
$u_n \to u$
uniformly in $\Omega$. 

Let now $\{F_n\}_{n>0}$ be a sequence of uniformly elliptic concave (or 
convex)
operators with ellipticity constants $\lambda$ and $\Lambda$ such
that $F_n \to F$ uniformly in compact sets of ${\cal S}_n \times
\Omega$ (${\cal S}_n$ is the set of symmetric matrices).  Suppose
in addition that $u_n$ satisfies
$$F_n(D^2u_n,x)=0\quad \mbox{in}\quad \Omega, \quad
u_n=0\quad\mbox{on}\quad\partial\Omega.$$
and that $\,u_n\,$ converges uniformly to $\,u$.
Then, $u \in C(\bar\Omega)$ is a solution to
$$F(D^2u,x)=0\quad \mbox{in}\quad \Omega, \quad
u=0\quad\mbox{on}\quad\partial\Omega .$$
\end{prop} 

\begin{Remark} Actually, the above proposition is proved in \cite{cab} in a 
more general case,
when $\{f_n\}_{n>0} \subset L^\infty(\Omega)$ and
$\{u_n\}_{n>0} \subset C(\bar\Omega)$ is a
sequence of viscosity solutions to \re{eqcom}. 

So, to prove Proposition \ref{pc} we need to use the following fact.
If $u \in  C(\bar\Omega) \cap W^{2,N}_{loc}(\Omega)$ is a sub-solution
(resp. super-solution) of ${\cal
M}^+(D^2u)=g$ with $g$ continuous, then $u$ is also viscosity sub-solution
(resp. super-solution) of the same equation, see \cite{CKS}.
We also use the regularity to prove that the limit of  $u_n$, $u$,  belongs to  
$\,W^{2,N}_{loc}(\Omega)$.
\end{Remark} 

\bigskip
Now we want to study the nonlinear bifurcation problem. We will first prove
the following. 

\begin{prop}\label{prop1} If $(\bar\mu, 0)$ is a bifurcation point of
problem \re{EVB},
then $\bar\mu $ is an eigenvalue of $\mmp$.
\end{prop} 

\noindent{\bf Proof.}
Since $(\bar\mu, 0)$ is a nonlinear bifurcation point, there is a sequence
$\{(\mu_n, u_n)\}_{n \in \NN}$
of nontrivial solutions of the problem \re{EVB} such that $\mu_n \to\bar\mu$
and $u_n \to 0$
in uniformly in $\Omega$.
Let us define
$$\hat u_n=\frac{u_n}{\|u_n\|_{\li}}$$
then $\hat u_n$ satisfies
$$-\mmp(D^2\hat u_n)= \mu_n \hat u +\frac{f(u_n, \mu_n)}{u_n}\hat
u_n\quad\mbox{in} \quad \Omega $$
So, the right-hand side of the equation is bounded.
Then by Proposition \ref{pc} we can extract a subsequence
such that $\hat u_n \to \hat u$. Clearly $\hat u$ is a solution to
\re{EV}.$\Box$ 

\bigskip
Before proving Theorem \ref{bteo}, we need some preliminaries  in order to compute
the Leray-Schauder degree of a related
function. 

\medskip
To start, let us recall some basic properties of the matrix operators
$\mmp$, whose proof follows directly from the equivalent definition for
$\mmp$:
$$\displaystyle \mmp(M)=\sup_{A \in {\cal A}}\mbox{tr}(AM),$$
for any symmetric matrix  $M$ (see \cite{cab}). Notice that the original
definition of
C. Pucci \cite{CPucci1} is of this type, but ${\cal A}$ is a different family
of symmetric matrices. 

\begin{lem}\label{basic}
Let $M$ and $N$ be two symmetric matrices. Then:
$$\mmp (M+N)\leq \mmp (M)+\mmp(N).$$
\end{lem} 

Next we recall a very well known fact about Pucci's operator,
namely that is the Alexandroff-Bakelman-Pucci estimate holds. The
proof can be found for example in \cite{cab}. 

\begin{teo} {\bf [ABP]}\label{teo4}
Let $\Omega$ be a bounded domain in $\RR^N$, such that\\
$\mbox{diam}(\Omega) \leq d$ and $f \in  L^N(\Omega)$.
Suppose that $u$ is continuous in $\overline\Omega$ and satisfies
$\mmm (D^2u) \leq f(x)$ in $\Omega$
and $u \leq 0$ on $\partial \Omega$.
Then,
$$\sup u^- \le C \|f^+\|_{L^N(\Omega)}.$$
Here $C=C(\mbox{meas}(\Omega),\lambda, \Lambda,N, d)$ is a constant and
meas$(\Omega)$ denotes Lebesgue
measure of $\Omega$.
\end{teo} 

The next corollary is a maximum principle for small domains, that was first
noted by Bakelman
and extensively used in \cite{berest} . 

\begin{cor}\label{cor1}  Let $\Omega$ be a bounded domain in $\RR^N$, such
that $\mbox{diam}(\Omega) \leq d$.
Suppose that  $u$ is continuous in $\overline\Omega$, satisfies $\mmm
(D^2u)+ c(x)u(x) \leq 0$ in $\Omega$, $u \geq 0$ on
$\partial\Omega$ and $c \in L^\infty(\Omega)$ with $c(x) \leq b$ a.e.. There
exists  $\delta=\delta(\lambda, \Lambda,N, d, b)$ such that
meas$(\Omega)< \delta$ implies $ u \geq 0$ in $\Omega$.
\end{cor} 

The proof is standard and uses in a crucial way Theorem \ref{teo4}. For
details see \cite{berest}.
Next corollary is crucial to prove that the
eigenvalue $\mu^-_1$ is isolated.
\medskip
\begin{cor}\label{cor2}
Let $\Omega_n$ be a sequence of domains such that meas$(\Omega_n)
\to 0$ as $n \to \infty$ and diam$(\Omega_n) \leq d$. If $(\mu_n
,u_n)$ is a positive solution to \re{EV} with $\Omega=\Omega_n$,
then $\mu_n \to \infty $ as $n \to \infty$.
\end{cor}
\noindent{\bf Proof:} Suppose by contradiction that there exists
$C>0$ such that $\mu_n< C$. Then $u_n$ satisfies the equation
$\mmp(D^2 u_n)+Cu_n \geq 0$. Since the measure of $\Omega_n$ is
small for $n$ large, we can use the previous corollary with $-u_n$
concluding that $-u_n \geq 0$, which is a contradiction. $\Box$ 

\begin{Remark}
1) In the sequel, we will vary the parameter $\lambda$ while
keeping $\Lambda$ fixed in the operator $\mmp$. We will denote the
half eigenvalues $\mu^+_1(\lambda)$, $\mu^-_1(\lambda)$ to make
explicit the dependence on the parameter $\lambda \in
(0,\Lambda]$. 

\bigskip
\noindent
2) From the characterization ii) of Proposition \ref{fet} it follows that if
$\lambda_1<\lambda_2$,
then  $\mu^+_1(\lambda_1) \leq \mu^+_1(\lambda_2)$ and
$\mu^-_1(\lambda_1) \geq \mu^-_1(\lambda_2)$. 

\end{Remark} 

\begin{lem}\label{contvp}
The two first half eigenvalues functions $\mu^+_1:(0,\Lambda]\to \RR$ and
$\mu^-_1:(0,\Lambda] \to \RR$,
are continuous on $\lambda$.
\end{lem} 

\noindent{\bf Proof:} Let $\{\lambda_j\}_{j \in \NN}$ be sequence in
$(0,\Lambda]$
converging to $\lambda \in (0,\Lambda]$.
We will show that
$$\displaystyle\lim_{j \to \infty }\mu^+_1(\lambda_j)=\mu^+_1(\lambda).$$
Since $\lambda_j \to \lambda$ there exists $\varepsilon>0$ such that
$\bar\lambda:=\lambda+\varepsilon \geq \lambda_j \geq
\lambda-\varepsilon=:\lambda^*>0$, for $j$ large.
 From the previous Remark we have
$$\mu^+_1(\lambda^*) \leq \mu^+_1(\lambda_j) \leq \mu^+_1(\bar\lambda),$$
for large $j$. Therefore, up to subsequences $\mu^+_1(\lambda_j) \to \mu$. 

Let $u_j$ be the corresponding eigenfunction for the eigenvalue
$\mu^+_1(\lambda_j)$. We can
suppose that $\|u_j\|_{\li}=1$,
then $u_j$ satisfies
$$\mmpj(D^2u_j)\geq -\mu^+_1(\bar\lambda) u_j \quad\mbox{and}\quad
\mmmj(D^2u_j)\leq -\mu^+_1(\lambda^*) u_j.$$
So by Proposition \ref{pc} up to subsequences, $u_j \to u$ uniformly in
$\Omega$.
Moreover, $(\mu,u)$ is a solution to \re{EV} and $\|u\|_{\li}=1$. 

Since $u_j$ is positive in $\Omega$,
we have that  $u_j$ is non-negative in $\Omega$ and by the strong maximum
principle, $u$ is positive in $\Omega$ .
Hence, by the uniqueness of the positive eigenfunction, Proposition \ref{fet}
i), $\mu=\mu^+_1(\lambda)$,
which ends the proof in this case.
The same proof holds in the case of $\mu^-_1$. $\Box$ 

\bigskip
The next Lemma proves that the first half-eigenvalue $\mu^-_1$ is isolated. 

\begin{lem}\label{lem3}
For every interval $[a,b] \subset (0,\Lambda)$ there exists a $\delta>0$ such
that for all
$\lambda \in [a,b]$ there is no eigenvalue of \re{EV} in
$(\mu^-_1(\lambda),\mu^-_1(\lambda)+\delta]$.
\end{lem} 

\noindent{\bf Proof:} Suppose that the Lemma is not true. Then, there are
sequences
$\{\lambda_j\}_{j \in \NN} \subset (0,\Lambda] $, $\{\mu_j\}_{j \in \NN}
\subset \RR^+$, and $\{u_j\}_{j \in \NN}
\subset C(\Omega)\setminus  \{0\}$
such that $\lambda_j \to \bar\lambda \in (0,\Lambda)$,  $\mu_j >
\mu^-_1(\lambda_j)$,
$\lim_{j \to \infty}(\mu_j-\mu^-_1(\lambda_j))=0$, and
$$-\mmpn(D^2 u_n)=\mu_n u_n.$$
Using Proposition \ref{pc} we have that, up to a subsequence, $u_n
\to u$ uniformly in $\Omega$ and $u$ is a solution of the problem
$$-\mmp(D^2 u)=\mu^-_1(\bar\lambda) u \quad\mbox{in}\quad\Omega.$$
Therefore by Proposition \ref{fet} i), $u$ is negative in $\Omega$. 

On the other hand, by i) of Proposition \ref{fet} $u_n$ changes sign
in $\Omega$, then there exists $\Omega_n$, a connected component
of $\{x \in \Omega\,|\,u_n(x) > 0\}$, with\\
$\mathrm{meas}(\Omega_n)>0$. Since $u_n \to u$,
$\mathrm{meas}(\Omega_n) \to 0$. Then by Corollary \ref{cor2}
$\mu(\Omega_n,\lambda^* ) \to \infty$, where $\lambda^*>0$ is such
that $\lambda_j >\lambda^*$. But $\mu_n=\mu^+_1(\Omega_n,\lambda_n
) \geq \mu(\Omega_n,\lambda^* )$, contradicting the fact that
$\mu_n$ converges to $\mu^-_1(\bar\lambda)$. $\Box$ 

\bigskip
Let us define
$$\mu_2(\lambda)=\inf\{\mu>\mu^-_1(\lambda)\,|\quad\mu \quad \mbox{is an
eigenvalue of \re{EVB}}\}$$
Then by the previous Lemma $\mu_2>\mu^-_1$. We notice that $\mu_2$
may be equal to $+ \infty$. Define now ${\cal L}^+_\lambda$ as the
inverse of $-\mmp$. It is well known that ${\cal L}^+_\lambda$ is
well defined in ${\cal C}:=\{u \in C(\bar\Omega)\,|\,u=0
\quad\mbox{on}\quad\partial\Omega\}$ (see for example \cite{caff})
and, by Proposition \ref{pc}, ${\cal L}^+_\lambda$ is compact.

Now we are in position to compute the Leray-Schauder degree and prove the
following proposition. 

\begin{prop}\label{prop2}
Let $r>0$, $\bar\lambda>0$, $\mu \in \RR$. Then
$$\mathop{\mathrm{deg}}_{{\cal C}}(I-\mu{\cal L}^+_{\bar \lambda}, B(0,r),
0)=\left\{\begin{array}{lll}
     1 & \mbox{if}\quad \mu<\mu^+_1(\bar \lambda)\\
     0 &  \mbox{if}\quad \mu^+_1(\bar \lambda)<\mu<\mu^-_1(\bar \lambda)\\
     -1 &\mbox{if}\quad \mu^-_1(\bar \lambda)< \mu <\mu_2(\bar \lambda),
\end{array}\right.$$
here ${\cal C}:=\{u \in C(\bar\Omega)\,|\,u=0
\quad\mbox{on}\quad\partial\Omega\}$.
\end{prop} 

\begin{Remark} Since ${\cal L}^+_\lambda$ is compact, the degree is well
defined if
$0 \not\in (I-\mu{\cal L}^+_{\bar \lambda})(\partial B(0,r)).$
\end{Remark} 

\noindent{\bf Proof of Proposition \ref{prop2}.}
We have that the degree
$$\mbox{deg}_{{\cal C}}(I-s\mu{\cal L}^+_{\bar \lambda}, B(0,r), 0)$$
is well defined for any $s \in [0,1]$ and $\mu<\mu^+_1(\bar\lambda)$,
since $\mmp$ does not have eigenvalues below $\,\mu_1^+$,  that is, $0 \not\in (I-s\mu{\cal 
L}^+_{\bar\lambda})(\partial B(0,r))$.
Using the invariance of the degree
under homotopy, we conclude that this degree is equal to 1, its
value at $s=0$. 

In the case $\mu^+_1(\bar \lambda)<\mu<\mu^-_1(\bar \lambda)$
we will use the following property of the degree to prove that the degree is
zero.
If $\mbox{deg}_{{\cal C}}(I-\mu{\cal L}^+_{\bar \lambda}, B(0,r), 0) \not
=0$, then
$(I-\mu{\cal L}^+_{\bar \lambda})(B(0,r))$ is a neighborhood of zero.
So we claim that if $\mu^+_1(\bar \lambda)<\mu<\mu^-_1(\bar \lambda)$,
then $(I-\mu{\cal L}^+_{\bar \lambda})(B(0,r))$ is not a neighborhood of
zero.
Suppose by contradiction that $(I-\mu{\cal L}^+_{\bar \lambda})(B(0,r))$ is
a neighborhood of zero.
Then for any smooth $h$ with $\|h\|_{\li}$ small,
there exists $u$ a solution to
$$u-\mu{\cal L}^+_{\bar \lambda}u=h$$
In particular, we can take $h$ to be a solution of
$$\mmpb(D^2h)=-\delta \quad\mbox{in}\quad \Omega\quad\mbox{and}\quad h=0
\quad\mbox{on}\quad\partial\Omega,$$
where  $\delta>0$ is small enough. 

Then, by Lemma \ref{basic} and the definition of ${\cal L}^+_{\bar
\lambda}$, it follows that $u$ satisfies
$$\mmpb(D^2u)+\mu u \leq -\delta\quad\mbox{in}\quad \Omega.$$
On the other hand, by Lemma \ref{lemh} and Remark \ref{rhopf} 3), 
there
exists $\varepsilon>0$, such that $\varepsilon(-\varphi)^-_1<u$,
and $\varepsilon(-\varphi^-_1)$ satisfies
$$\mmpb(D^2 \varepsilon (-\varphi^-_1))+\mu \varepsilon (-\varphi^-_1) \geq
 -\delta\quad\mbox{in}\quad \Omega.$$
Then using Perron's method we find a positive solution $w$ to
$$\mmpb(D^2 w)+\mu w =-\delta\quad\mbox{in}\quad \Omega,\quad w=0
\quad\mbox{on}\quad\partial\Omega.$$
This leads to a contradiction with Theorem \ref{fet1} and with
the characterization for the first eigenvalue \re{rem12}
1). So, $\mbox{deg}_{{\cal C}}(I-\mu{\cal L}^+_{\bar \lambda},
B(0,r), 0)=0$ for $\mu^+_1(\bar \lambda)<\mu<\mu^-_1(\bar
\lambda)$. 

Finally, suppose that $\mu^-_1(\bar \lambda)< \mu <\mu_2(\bar \lambda)$. The
continuity of $\mu^-_1(\cdot) $
and Lemma \ref{lem3} imply the existence of a continuous function
$\nu:(0,\Lambda] \to \RR$ such that $\mu^-_1(\lambda)< \nu(\lambda)
<\mu_2(\lambda)$ for all $\lambda \in (0,\Lambda]$
and $\nu(\bar \lambda)=\mu$. 

The result will follow by showing that the well-defined, integer-valued
function
$$d(\lambda)=\mbox{deg}_{{\cal C}}(I-\nu(\lambda){\cal L}^+_{\bar \lambda},
B(0,r), 0)$$
is constant in $[\bar \lambda,\Lambda]$.
This follows by the invariance of the Leray-Schauder degree under a compact
homotopy.
Recall that $d(\Lambda)=-1$, hence the proposition follows. $\Box$ 

{\bf Proof of Theorem \ref{bteo}.}  Let us set
$$H_\mu(u)={\cal L}^+_\lambda(\mu u+f(\mu, u)).$$ 

Suppose that $(\mu^+_1,0)$ is not a bifurcation point of problem \re{EVB}.
Then there exist $\varepsilon$, $\delta_0>0$ such that for all
$|\mu-\mu^+_1|\leq \varepsilon$
and $\delta < \delta_0$ there is no nontrivial solution of the equation
$$u-H_\mu(u)=0$$
with $\|u\|=\delta$. From the invariance of the degree under compact
homotopy we obtain that
\beq \label{grad1}
\mbox{deg}_{{\cal C}}(I-H_\mu, B(0,\delta), 0) \equiv \mbox{constant}\quad
\mbox{for} \quad \mu
\in [\mu^+_1-\varepsilon,\mu^+_1+\varepsilon]
\eeq
By taking  $\varepsilon$ smaller if necessary, we can assume that
$\mu^+_1+\varepsilon < \mu^-_1$.
Fix now $\mu \in (\mu^+_1, \mu^+_1+\varepsilon]$.
It is easy to see that if we choose $\delta$ sufficiently small, then the
equation
$$u-{\cal L}^+_\lambda(\mu u+sf(\mu, u))=0$$
has no solution $u$ with $\|u\|=\delta$ for every $s \in [0,1]$.
Indeed, assuming the contrary and reasoning as in the proof of Proposition
\ref{prop1}, we would find that
$\mu$ is an eigenvalue of \re{EV}. From the invariance of the degree under
homotopies and Proposition \ref{prop2} we obtain
\beq\label{grad2}
\mbox{deg}_{{\cal C}}(I-H_\mu, B(0,\delta), 0)=\mbox{deg}_{{\cal
C}}(I-\mu{\cal L}^+_\lambda , B(0,\delta), 0)=0
\eeq
Similarly, for $\mu \in [\mu^+_1-\varepsilon, \mu^+_1)$ we find that
\beq\label{grad3}
\mbox{deg}_{{\cal C}}(I-H_\mu, B(0,\delta), 0)=1
\eeq
Equalities \re{grad2} and \re{grad3} contradict \re{grad1} and hence
$(\mu^+_1, 0)$ is a bifurcation point for the
problem \re{EVB}. Let define $u_\mu$ a solution to \re{EVB} for
$\mu>\mu^+_1$, with
$\|u_\mu\|_\infty \to 0$ as $\mu \to \mu^+_1$. Using the same argument of
Proposition \ref{prop2},
$$u_\mu /\|u_\mu\|_\infty \to \varphi^+_1 \quad \mbox{as} \quad \mu \to
\mu^+_1.$$
This shows that $u_\mu$ is positive for $\mu$ close to $\mu^+_1$. 

The rest of the proof is entirely similar to that of the Rabinowitz's Global
Bifurcation Theorem,
see \cite{rabinowitz2}, \cite{rabinowitz} or \cite{rabinowitz1}, so we omit it here. $\Box$ 




\section{``Spectrum'' in the Radial Case and Nonlinear Bifurcation from all
``eigenvalues''} 

\medskip 

Let us first recall that the value of the Pucci's operator applied
to a radially symmetric function can be computed explicitly;
namely if $u(x)=\varphi(|x|)$ one has
$$
D^2 u  (x)=\frac{\varphi '(|x|)}{|x|}I+\left[\frac{\varphi ''(|x|)}{|x|^2}-
\frac{\varphi '(|x|)}{|x|^3}\right] x\otimes x ,
$$
where $I$ is the $N\times N$ identity matrix and
$x \otimes x$ is the matrix whose entries are $x_ix_j$. Then the eigenvalues
of
$D^2u$ are $\varphi ''(|x|)$, which is simple, and
${\varphi '(|x|)}/{|x|}$, which has multiplicity $N-1$. 

In view of this, we can give a more explicit definition of Pucci's operator.
In the case of $\mmp$ we define the functions
\begin{eqnarray*}
M(s)=
\left\{\begin{array}{ll}
{s}/{\Lambda} \quad s>0,\\
{s}/{\lambda} \quad s\le 0,
\end{array}\right.
\quad\mbox{ and }\quad
m(s)=
\left\{\begin{array}{ll}
\Lambda s \quad s>0,\\
\lambda s \quad s\le 0.
\end{array}\right.
\end{eqnarray*}
Then, we see that $u$ satisfies  \equ{EV} with $\,\Omega=B_1\,$  and is radially symmetric if and
only if $u(x)=v(|x|)$, $r=|x|$
satisfies
\be\label{eqv}
v''= M(-\frac{(N-1)}{r} m(v')-\mu v),
\ee
\be\label{eqvbc}
v'(0)=0,\quad v(1)=0.
\ee 

Next we briefly study the existence, uniqueness, global existence,
and oscillation of the solutions to the related initial value
problem \be\label{eqv1} w''= M(-\frac{(N-1)}{r} m(w')- w), \ee
\be\label{eqvbc1} w'(0)=0,\quad w(0)=1. \ee Then we will come back
to \re{eqv}, \re{eqvbc} and to the proof of Theorem \ref{teovpr}.
First using a standard Schauder fixed point argument as used by Ni
and Nussbaum in \cite{Ninuss}, we can prove the existence of $w
\in C^2$ solution to
$$\{w'r^{N-1}\}'=-r^{N-1}\frac{w}{\lambda} , \quad w'(0)=0,\quad w(0)=1. $$
Moreover, this solution is unique and  for $r$ small, $w'(r)$ and $w''(r)$ are negative.
Then, for some $\delta>0$,  $w$ satisfies
$$w''= M(-\frac{(N-1)}{r} m(w')-  w),\quad\mbox{in}\quad (0, \delta]$$
Next we consider \equ{eqv1} with initial values $w(\delta)$ and
$w'(\delta)$ at $r=\delta$. From the standard theory of ordinary
differential equations we find a unique $C^2$-solution of this
problem for $r\in[\delta,a)$, for $a>\delta$. Using Gronwall's
inequality we can extend the local solution to $[0,+\infty)$. 

In the following Lemma we will show that the solution $w$ is oscillatory. 

\begin{lem}The unique solution $w$ to \re{eqv1} \re{eqvbc1}, $w$, is
oscillatory, that is,
given any $r>0$, there is a $\tau>r$ such that $w(\tau)=0$.
\end{lem}
The proof uses standard arguments of oscillation theory for ordinary
differential equation. 

{\bf Proof.} Suppose that $w$ is not oscillatory, that is, for some $r_0$ ,
$w$ does not
vanish on $(r_0, \infty)$.
Assume that  $w>0$ in $(r_0, \infty)$. Let $\phi$ be a solution to
\re{eqv1}, \re{eqvbc1} with
$\lambda=\Lambda$, then it is known that $\phi$ is oscillatory.
So we can take $r_0<r_1,\, r_2$ such that $\phi(r)>0$ if $r \in (r_1,r_2)$
and $\phi(r_1)=\phi (r_2)=0$.
We have that $w$ and $\phi$
satisfy
$$\{ w'r^{N-1}\}'\leq -r^{N-1}\frac{w}{\lambda},$$
$$\{ \phi'r^{N-1}\}'= -r^{N-1}\frac{\phi}{\lambda}.$$
If we multiply the first equation by $\phi$ and the second by $w$, subtract
them and then integrate, we get
$$r_1^{N-1}\phi'(r_1) w(r_1)-r_2^{N-1}\phi'(r_2) w(r_2) \leq 0,$$
getting a contradiction. 

Suppose now that $w<0$ in $(r_0, \infty)$.
In that case we claim that  $w'> 0$ in $(r_0, \infty)$, taking if necessary a larger
$r_0$. If there exists a $r^*$ such that $w' (r^*)=0$, then using the
equation
we have that $w'> 0$ in $(r^*, \infty)$. So we only need to discard the case
$w'< 0$ in $(r_0, \infty)$. In that case $w$ satisfies
$$\{ w'r^{\tilde N^+-1}\}'= -r^{\tilde N^+-1}\frac{w}{\lambda} \quad
\mbox{in}\quad (r_0, \infty)$$
where $\tilde N^+=(\lambda(N-1))/\Lambda+1$.
Let denote by $g(r)=\{w'r^{\tilde N^+-1}\}$ we have that $g$ is monotone,
then
there exists a finite $c_1<0$ such that $\lim_{r\to\infty}g(r)=c_1$. 

On the other hand, since $w'< 0$, there exists $c_2 \in [-\infty, 0)$
such that $\lim_{r\to\infty}w(r)=c_2$, then from the equation satisfied by
$w$, we get that $$\displaystyle \lim_{r\to\infty} g'(r)=+\infty.$$
That is a contradiction with
$\displaystyle  \lim_{r\to\infty}g(r)=c_1$. 

Define now
$$b(r)=r^{\tilde N^--1}\frac{w'(r)}{w(r)},\quad r \in (r_0, \infty),$$
here $\tilde N^-=(\Lambda(N-1))/\lambda+1$. 

Then we claim that $b$ satisfies,
\be\label{os}
b'+\frac{b^2}{r^{\tilde N^--1}}+\frac{r^{\tilde N^--1}}{\Lambda} \leq 0.
\ee
If $w''>0$ then  $b$ satisfies
$$
b'+\frac{b^2}{r^{\tilde N^--1}}+\frac{r^{\tilde N^--1}}{\lambda}=0
$$
Since $\frac{1}{\Lambda} \leq \frac{1}{\lambda}$, the claim follows in this
case.
If $w''<0$ then  $b$ satisfies
$$b'+\frac{b^2}{r^{\tilde N^--1}}+\frac{r^{\tilde N^--1}}{\Lambda} = (\tilde
N^--N)b.$$
Finally, since $\tilde N^--N \geq 0$ and $b<0$, the claim follows also in
this second case. 

Integrating \re{os} from $r_0$ to $t>r_0$ we get
\be \label{inte}
b(t)-b(r_0)+\frac{t^{\tilde N^-}}{\tilde N^-\Lambda}-\frac{r_0^{\tilde
N^-}}{\tilde N^-\Lambda}
+\int^t_{r_0} \frac{b^2}{r^{\tilde N^--1}} \leq 0.
\ee
In particular we have
$$-b(t) \geq Ct^{\tilde N^-}.$$
For some $C>0$ and $t$ large.
Define now
$$k(t)=\int^t_{r_0} \frac{b^2}{r^{\tilde N^--1}}.$$
Then, by the previous fact, we have
\be \label{asi1}
k(t)\geq ct^{\tilde N^-+2 }\quad \mbox{for } t\mbox{ and some } c>0.
\ee
On the other hand from \re{inte} and $b<0$ we get
$$k(t)<-w(t),$$
or
$$k(t)<k'(t)t^{\tilde N^-+1 }, \quad \mbox{for } t\mbox{ large. }$$
The latter inequality implies
\be \label{kk}
C(\frac{1}{k(t)}-\frac{1}{k(s)}) \geq
\frac{1}{t^{N^--2}}-\frac{1}{s^{N^--2}}.
\ee
for some $C>0$ and $t,s$ large with $t<s$. Letting $s \to \infty$ and
noting that $k(s) \to +\infty$, we find
\be \label{asi2}
k(t) \leq A t^{N^--2}.
\ee 

However \re{asi1} and \re{asi2} are not compatible. This contradiction shows
that $w$ must be oscillatory. $\Box$ 

\bigskip
Notice that the same proof holds when the initial conditions to the problem
\re{eqv1}
are $w(0)=-1$, $w'(0)=0$. 

\bigskip
With these preliminaries we are now ready to prove Theorem \ref{teovpr}. 

\bigskip
{\bf Proof of Theorem \ref{teovpr}.} Let denote $w^\nu$ the above solutions
of \re{eqv1}
with initial conditions $w^\nu(0)=\pm 1$ (here and in the rest of the proof
$\nu \in \{+,\,-\}$).
 From the previous lemma, $w^\nu$ has infinitely many zeros:
$$0<\beta^\nu_1<\beta^\nu_2<\cdot\cdot\cdot<\beta^\nu_k<\cdot\cdot\cdot.$$
A standard Hopf type argument shows that they are all simple. Next
we define $\mu^\nu_k=(\beta^\nu_k)^2$ of Theorem \ref{teovpr}. Clearly
$\mu= \mu^\nu_k$ is an eigenvalue of \re{EV}, with $w^\nu(\beta^\nu_k
\,\cdot)$, $r \in [0,1]$, being the corresponding eigenfunction with
$k-1$ zeros in $(0,1)$. We claim that there is no radial eigenvalue of
\re{EV} other than these $\mu^\nu_k$'s. 

Let $\mu$ be an eigenvalue of \re{EV}. Clearly $\mu>0$. Let $z(r)$ be the
corresponding eigenfunction and suppose that 
$z(0)>0$, the uniqueness of solution to \re{eqv1} implies that
$z(r)=z(0)\, w^+(\mu^{1/2}r)$.
Moreover, since $z(1)=0$, $\mu=(\beta^+_k)^2$ for some $k \in \NN$, and
$z=z(0) \,w^+$. The same holds for $z(0)<0$. $\Box$ 

\bigskip
Below we will exhibit some properties of the eigenvalues distribution. 

\begin{lem}\label{lemorden}
For $k \in \NN,\,k>1$ we have
$\mu^-_k < \mu^+_{k+1}$
and
$\mu^+_k < \mu^-_{k+1}$.
\end{lem} 

\noindent
{\bf Proof.}
We will prove the lemma in terms of the functions $w^+$  and
$w^-$ defined above. 

We claim that if $w^+$ has to change sign between two
consecutive zeros of $w^-$, if $w^+$ has the same sign of $w^-$. 
Notice that this is weaker then the usual Sturm's comparison result, since 
there
is a additinal sign restriction. 

Suppose first by contradiction that $w^-(r_1)=w^-(r_2)=0$, $w^-(r)>0$ for 
all $r
\in (r_1,r_2)$
and $w^+(r)>0$ for all $r \in [r_1,r_2]$.
Let $r_3<r_1<r_2<r_4$ be the next zeros of $w^+$, that is,
$w^+(r_3)=w^+(r_4)=0$, $w^+(r)>0$ for all $r \in (r_3,r_4)$.
Then, the first half-eigenvalue in $A_1:=\{r_1<|x|<r_2\}$ is $\mu^+(A_1)=1$
and
first half-eigenvalue in $A_2:=\{r_3<|x|<r_4\}$ is $\mu^+(A_2)=1$.
Define now $u(r)=w^+(\beta r)$, with $\beta>1$ such that $r_4/\beta > r_2$.
So,
$u$ is a positive eigenfunction in
$A_3:=\{\frac{r_3}{\beta}<|x|<\frac{r_4}{\beta}\}$ with
eigenvalue $\mu^+(A_3)=\beta^2$. But $A_1 \subset A_3$, therefore
$\mu^+(A_1)=1\geq \mu^+(A_3)=\beta^2$
getting a contradiction. The same kind of argument can be used in the case 
when
$w^-$ negative in $(r_1, r_2)$ and $w^+$  negative in $[r_1,r_2]$.
Hence, the claim follows. 

In the two cases above we can invert the role of $w^-$
and $w^+$. 

As a consequence of the previous facts, the lemma follows by
examining the distribution of zeroes of $w^+$ and $w^-$.$\Box$ 

\begin{Remark}\label{rzero}
1) The above lemma implies that in the case
$\beta^+_k <\beta^{-}_k$,\\
$w^+(r)\,w^-(r)>0$ for all $r \in (\beta^+_k ,\beta^{-}_k)$.
The same holds true in the case $\beta^+_k >\beta^{-}_k$.
\end{Remark} 

\begin{lem}
The gap between the two first half-eigenvalues is larger than that betwenn 
the second ones:
$$\frac{\mu^-_1}{\mu^+_1} \geq \frac{\mu^-_2}{\mu^+_2}.$$
\end{lem} 

\noindent
{\bf Proof.} Let $\varphi^+_2$ and $\varphi^-_2$ the radial eigenfunctions 
of $\mmp$ in $B_1$, with corresponding eigenvalues
$\mu^+_2$ and $\mu^-_2$. Define $r^+$ (resp. $r-$) as the first zeros of  
$\varphi^+_2$ (resp. $\varphi^-_2$).
We claim that $r^- \geq r^+$. Suppose by contradiction that $r^- < r^+$. 
Define now
$A^+=\{x\,|\,r^+<|x|<1\}$ and $A^-=\{x\,|\,r^-<|x|<1\}$, then $A^+ \subset 
A^-$.
Using the monotonicity with respect the domain of the first half-eigenvalues 
and Proposition \ref{fet} ii) we get
$$\mu^-_1(A^+)=\mu^+_2 \geq \mu^+_1(A^+) \geq \mu^+_1(A^-)= \mu^-_2 .$$
On the other hand $B_{r^-} \subset B_{r^+}$, thus by the same kind of 
argument
$$\mu^-_1(B_{r^-})=\mu^-_2 > \mu^-_1(B_{r^+}) \geq \mu^+_1(B_{r^+}) = 
\mu^+_2 .$$
Hence, we get a contradiction . So, the claim follows.
Making a rescaling argument, so as in the proof of Theorem \ref{teovpr}, it 
follows that
$$(r^+)^2  \mu^+_2 = \mu^+_1 \quad\mbox{and}\quad (r^-)^2  \mu^-_2 = 
\mu^-_1,$$
which ends the proof. $\Box$

\bigskip
Next, we prove some preliminary results to prepare the proof of Theorem 
\ref{teobr}.

\begin{lem}\label{lemne} Assume that $\mu^+_k\not = \mu^-_k$ and that there 
exists $r_0 \in (0,1)$ such that
$\phi_\pm(r)>0$ for all $r \in (r_0,1]$.
Then, there exists a continuous function $g$ such that there is no solution 
to the problem
\be\label{uf}
u''= M(-\frac{(N-1)}{r} m(u')-\mu u+g)\quad\mbox{in}\quad [0,r_0],
\ee
and
\be\label{uf1}
u'' \geq M(-\frac{(N-1)}{r} m(u')-\mu u+g)\quad\mbox{in}\quad (r_0,1],
\ee
\be\label{cbf}
u'(0)=0,\quad u(1)=0.
\ee
for $\mu $ between $\mu^+_k$ and $\mu^-_k$.
\end{lem}
\begin{Remark}1) Some ideas of the proof are in the book of P. Dr\'abek 
\cite{dra}. 

\bigskip
\noindent
2) There is a similar non-existence result in the case when there exists 
$r_0 \in (0,1)$ such that
$\phi_\pm(r)< 0$ for all $r \in (r_0,1]$ , replacing \re{uf1} by
\be\label{uf2}
u'' \leq M(-\frac{(N-1)}{r} m(u')-\mu u+g) \quad\mbox{in}\quad (r_0,1],
\ee
in the previous lemma. 

\bigskip
\noindent
3) Let us denote by $\phi_+$ and $\phi_-$ the solutions of \re{uf} with $\,r_0=1$ 
and $g=0$
 and respective initial conditions
$u'(0)=0,$ $u(0)=1$ and $u'(0)=0,$ $u(0)=-1$.
Let us suppose that $\mu$ is between $\mu^+_k$ and $\mu^-_k$, then by
Remark \ref{rzero} we deduce that  $\phi_+(1)\phi_-(1)>0$.

\end{Remark}

\noindent
{\bf Proof.} Consider then the particular case
$$\phi_\pm(r)> 0, \quad \phi'_\pm(r)\leq 0 \quad\mbox{for all}\quad r\in 
(r_0,1].$$
All other cases can be treated similarly. 

Let $g:[0,1] \to \RR$ be a continuous function such that $g(r)=0$ for all $r 
\in  [0,\,r_0]$
and  $g(r)>0$ for all $r \in (r_0,\,1]$. 

For $\alpha \in \RR$, let $\varphi_\alpha$ be the solution to \re{uf}, 
\re{uf1} and \re{cbf} with
$\varphi_\alpha(0)=\alpha$. 
For $\alpha > 0$,
we have
$$\varphi_\alpha(r)=\alpha \phi_+(r) \quad\mbox{for all}\quad r\in 
[0,r_0],$$
since uniqueness holds when $g=0$.
Put $r_1=\inf\{r\in (r_0, 1);\,\varphi_\alpha(r)= 0  \}$.
The interval $(r_0, r_1)$ contains a point $\tau_1$
such that
$$\displaystyle\left[\frac{\varphi_\alpha}{\phi_+}\right]'(\tau_1)<0$$
If this is not the case,
$$\displaystyle \frac{\varphi_\alpha(\tau)}{\phi_+(\tau)} \geq
\frac{\varphi_\alpha(r_0)}{\phi_+(r_0)}=\alpha>0, \quad
\tau \in (r_0, r_1 ),$$
which is impossible. So, we obtain
$$(\varphi_\alpha'\phi_+-\varphi_\alpha\phi_+')(\tau_1)< 0.$$
Define
$$G_i(r)=r^{\tilde{N}_i-1}(\varphi_\alpha'\phi_+-\varphi_\alpha\phi_+'),\quad
i=1,2,$$
where $\tilde{N}_1=N$ and $\tilde{N}_2=\tilde{N}^+$. 

Now we claim that there exists $\tau_2$, $r_0\leq \tau_2<\tau_1$  such that
$$\varphi_\alpha'(r)<0 \quad\mbox{for all}\quad r\in
(\tau_2,\tau_1)\quad\mbox{and}\quad G_i(\tau_2)\geq 0,\quad i=1,2.$$
If $\phi_\alpha'(r)<0$ for all  $r \in (t_0,\tau_1)$, since $G_i(r_0)=0$,
we conclude in this case by taking $\tau_2=r_0$.
If not, we define
$\tau_2=\sup\{\tau \in [r_0,\tau_1), \,|\,  \varphi_\alpha'(\tau)=0 \}.$
Notice that $\tau_2<\tau_1 $ and $\varphi_\alpha'(\tau_1)<0$, so
$\phi_\alpha'(r)<0$ for all $ r\in (\tau_2,\tau_1)$.
By the definition of $\tau_2$, $\varphi_\alpha'(\tau_2)=0$. Thus,
$G_i(\tau_2)>0$ and the claim follows. 
From the equation satisfied by $\phi_+$ we get
\beq\label{eqtau1}
\{r^{N-1}\phi'_+\}' \leq\frac{r^{N-1}}{\lambda}[ -\mu \phi_+
]\quad\mbox{in}\quad (\tau_2,\tau_1),
\eeq
and
\beq\label{eqtau2}
\{r^{\tilde N^+-1}\phi'_+ \}' \leq \frac{r^{\tilde N^+-1}}{\Lambda}
[ -\mu \phi_+ ]\quad\mbox{in}\quad (\tau_2,\tau_1).
\eeq
Since $\varphi_\alpha $ is positive in $(\tau_2,\tau_1)$, we obtain
$\displaystyle G'_1(r)\geq \frac{r^{N-1}}{\lambda}g(r)\phi_+(r)>0$, if
$\varphi_\alpha ''(r) <0$ and
$\displaystyle G'_2(r)\geq \frac{r^{\tilde N^+-1}}{\Lambda}
g(r)\phi_+(r)>0$, if $\varphi_\alpha ''(r)\geq 0$ for all $r \in
(\tau_2,\tau_1)$. 

\smallskip
The interval $(\tau_2, \tau_1)$ can be splitted
in subintervals $(s, t)$ such that
$G_i(s)-G_i(t)=\int^s_t G_i'(\tau)d\tau>0$, where $i$ is well chosen.
Using that if $G_i(t)<0$, then  $G_j(t)<0$ for $i \not= j$, we get a
contradiction. 

\smallskip
For $\alpha = 0$, $\varphi(r)=0$, $r \in [0, r_0]$. Then, we find an 
appropiate
interval to argue as
in the  above case.
For $\alpha < 0$ we have $\varphi_\alpha(r)=|\alpha |\phi_1$ for all
$r \in [0, r_0]$
and the proof is quite analogous as for  $\alpha > 0$.
All the above shows that there is no solution for
\re{uf}, \re{uf1} and \re{cbf}. $\Box$

\begin{prop}\label{prop3}
Let $r>0$, $\bar\lambda>0$, $\mu \in \RR$. Then
$$\mbox{deg}_{{\cal C}}(I-\mu{\cal L}^+_{\bar \lambda}, B(0,r),
0)=\left\{\begin{array}{lll}
     1 & \mbox{if}\quad \mu < \mu^+_1(\bar \lambda)\\
     0 &  \mbox{if}\quad \mu^+_k(\bar \lambda)<\mu<\mu^-_{k}(\bar
\lambda)\\
        &  \quad\mbox{or}\quad\mu^-_k (\bar \lambda)<\mu<\mu^+_k(\bar
\lambda)\\
     (-1)^k &\mbox{if}\quad \mu^+_k(\bar \lambda)< \mu <\mu^-_{k+1}(\bar
\lambda)\\
          &  \quad\mbox{or}\quad\mu^-_k(\bar \lambda)<\mu<\mu^+_{k+1}(\bar
\lambda),
\end{array}\right.$$
here ${\cal C}:=\{u \in C([0,1])\,|\,u(1)=0\,,\, u'(0)=0\}$.
\end{prop} 

\begin{Remark}
1) For $k \in \NN,\,k>1$, we do not expect that in general
$$\mu^+_k \leq \mu^-_k ,$$
but this is an open problem. 

\smallskip
\noindent
2) If $\mu^+_k = \mu^-_k $, the case
$\mbox{deg}_{{\cal C}}(I-\mu{\cal L}^+_{\bar \lambda}, B(0,r), 0)=0$
is not present in Proposition \ref{prop3}.
\end{Remark} 

\noindent
{\bf Proof.}
Assume first that $\mu^+_k(\bar \lambda)< \mu <\mu^-_{k+1}(\bar \lambda)$ or
$\mu^-_k(\bar \lambda)<\mu<\mu^+_{k+1}(\bar \lambda)$.
The arguments used in the proof of Lemma \ref{contvp} imply  
$\mu^\pm_j(\lambda)$ is
a continuous function
of $\lambda$. Using Lemma \ref{lemorden} we find a continuous function
$\nu:(0, \Lambda] \to \RR$
such that
$\max\{\mu^+_{k}(\lambda),\,\mu^-_{k}(\lambda)\}<\nu(\lambda)<\min\{\mu^+_{k
+1}(\lambda),\,\mu^-_{k+1}(\lambda)\}$ and $\nu(\bar \lambda)=\mu$. The
invariance of the Leray-Schauder's degree under compact homotopies
implies
$$d(\lambda)=\mbox{deg}_{{\cal C}}(I-\nu(\lambda){\cal L}^+_{\bar \lambda},
B(0,r), 0)=\mbox{constant},$$
for $\lambda \in (0,\Lambda]$. In particular $d(\bar
\lambda)=d(\Lambda)=(-1)^k$ and the result follows.
The case $\mu < \mu^+_1(\bar \lambda)$ is proved in Proposition \ref{prop2}.
In the case $\mu^+_k(\bar \lambda)<\mu<\mu^-_{k}(\bar \lambda)$ or $\mu^-_k
(\bar \lambda)<\mu<\mu^+_k(\bar \lambda)$
we will prove, as in Proposition \ref{prop2}, that
$(I-\mu{\cal L}^+_{\bar \lambda})(B(0,r))$ is not a neighborhood of zero. 

Suppose by contradiction that $(I-\mu{\cal L}^+_{\bar \lambda})(B(0,r))$ is
a neighborhood of zero.
Then, for any smooth $h$ with $\|h\|_{\lu}$ small,
there exists a solution $u$ to
$$u-\mu{\cal L}^+_{\bar \lambda}u=h$$
In particular we can take $h$ being a solution to
$$\mmpb(D^2h)=\psi \quad\mbox{in}\quad \Omega\quad\mbox{and}\quad h=0
\quad\mbox{on}\quad\partial\Omega,$$
where  $\|\psi\|_{\lu}>0$ is small enough.
Then, by Lemma \ref{basic} and the definition of ${\cal L}^+_{\bar
\lambda}$, it follows that $u$ satisfies
$$\mmpb(D^2u)+\mu u \leq \psi \quad\mbox{in}\quad \Omega.$$
Taking $\psi=-g$ (resp. $\psi=g$), where $g$ is a function of the type used 
in Lemma \ref{lemne},
we will get that $-u$ (resp. $u$) satisfies \re{uf}, \re{uf1} (resp.\re{uf2}) and  \re{cbf}.
Thus, we get a contradiction with lemma \ref{lemne} or Remark 3.2 2).
So, $\mbox{deg}_{{\cal C}}(I-\mu{\cal L}^+_{\bar \lambda}, B(0,r), 0)=0$ in
this case, and the proof is finished.$\Box$ 

\bigskip
\noindent
{\bf Proof of Theorem \ref{teobr}:} Using the same argument as in Theorem
\ref{bteo},
we obtain the existence of a ``half-component'' ${\cal C}^+_k \subset \RR 
\times
C([0,1])$ of radially symmetric
solutions to \re{EVB}, whose
closure $\bar{\cal C}^+_k$ contains $(\mu^+_k, 0)$ and is either unbounded
or contains a point
$(\mu^\pm_j, 0)$, with $j\not =k$ in the case of $\mu^+_j$. 

Let us first prove that if $(\mu,v) \in {\cal C}^+_k $, it implies that $v$
is positive at the origin
and possesses $k-1$ zeros in $(0,1)$. Arguing as in the proof of Theorem
\ref{bteo}, we find a
neighborhood ${\cal N}$ of $(\mu^+_k, 0)$ such that ${\cal N} \cap {\cal
C}^+_k \subset S^+_k$. 

Moreover, if $u \in C^1[0,1]$ is a solution to
\be\label{ufu}
u''= M(-\frac{(N-1)}{r} m(u')-\mu u+f(u, \mu))\quad \mbox{in}\quad(0,1).
\ee
and there exists $r_0 \in [0,1]$ such that $u(r_0)=u'(r_0)=0$, then $u
\equiv 0$. 

Using this fact we can extend the previous local properties of ${\cal
C}^+_k$ to all
of it. Hence, ${\cal C}^+_k$ must be unbounded. $\Box$

\end{document}